\documentclass[reqno,centertags, 12pt]{amsart}
\usepackage{amsmath,amsthm,amscd,amssymb,mathtools}
\usepackage{esint}  %from Maxim for circular integral
\usepackage{latexsym}
\usepackage{graphicx}
\usepackage{enumitem}
\usepackage[mathscr]{eucal}
\newcommand{\bbC}{{\mathbb{C}}}

\newcommand{\bbQ}{{\mathbb{Q}}}
\newcommand{\bbR}{{\mathbb{R}}}

\newcommand{\bbZ}{{\mathbb{Z}}}

\newcommand{\calH}{{\mathcal H}}

\newcommand{\calR}{{\mathcal R}}

\newcommand{\lb}{\label}

\newcommand{\wti}{\widetilde  }

\newcommand{\tr}{\text{\rm{Tr}}}

\newcommand{\spec}{\text{\rm{spec}}}

\newcommand{\bi}{\bibitem}

\newcommand{\beq}{\begin{equation}}
\newcommand{\eeq}{\end{equation}}
\newcommand{\ba}{\begin{align}}
\newcommand{\ea}{\end{align}}

\let\det=\undefined\DeclareMathOperator{\det}{det}

\newcounter{smalllist}

\newcommand{\comm}[1]{}
\DeclareMathOperator{\Ima}{Im}

\allowdisplaybreaks
\numberwithin{equation}{section}
\newtheorem{theorem}{Theorem}[section]

\theoremstyle{definition}
\newtheorem*{remark}{Remark}
\newtheorem*{remarks}{Remarks}

\newcommand{\jap}[1]{\langle #1 \rangle}
\newcommand{\bigjap}[1]{\left\langle #1 \right\rangle}
\newcommand{\Norm}[1]{\lVert#1\rVert}

\begin{document}

\title[Fritz and Me]{Fritz and Me---Perspectives in Inverse Spectral Theory}
\author[B.~Simon]{Barry Simon$^{1}$\\Departments of Mathematics and Physics\\ Mathematics 253-37\\ California Institute of Technology\\ Pasadena, CA 91125, USA.}
\thanks{$^1$ E-mail: bsimon@caltech.edu.}
\thanks{To Fritz Gesztesy on his $70^{th}$ birthday}

\date{\today}

\begin{abstract} This is a pedagogic introduction to certain aspects of inverse spectral theory for Schr\"{o}dinger operators and Jacobi matrices that revolves around my joint work with Fritz Gesztesy whose $70^{th}$ birthday we are honoring.
\end{abstract}

\maketitle

%%%%%%%%%%%%%%%%%%%%%%%%%%%%%%%%%%%%%%%%%%%%%%%%%%%%%%%%%%%%%%
\section{Introduction} \lb{s1}
%%%%%%%%%%%%%%%%%%%%%%%%%%%%%%%%%%%%%%%%%%%%%%%%%%%%%%%%%%%%%%

I have written 31 papers \cite{Fritz01, Fritz02, Fritz03, Fritz04, Fritz05, Fritz06, Fritz07, Fritz08, Fritz09, Fritz10, Fritz11, Fritz12, Fritz13, Fritz14, Fritz15, Fritz16, Fritz17, Fritz18, Fritz19, Fritz20, Fritz21, Fritz22, Fritz23, Fritz24, Fritz25, Fritz26, Fritz27, Fritz28, Fritz29, Fritz30, Fritz31} with Fritz Gesztesy (I still remember how delighted he was when he surpassed Elliott Lieb as my most frequent coauthor).  I will not attempt to review all of them here, but if there is a common thread in much of our work, it is inverse problems, how to use spectral data to determine the fundamental coefficients that specify the system - in the context of Schr\"{o}dinger  operators, most often the location of eigenvalues and their weights plus scattering data back to the potential, $V$. More generally, the basic issue is going from the Weyl-Titchmarsh $m$-function to $V$. A second, closely related issue that recurs in our work is spectral consequences of rank one perturbations.

Basic to the subject of rank one perturbations is the Aronszajn-Krein formula (after \cite{Ar1, Ar2, Ar3, Krein}).  We consider a bounded self-adjoint operator, $A$, on a Hilbert space, $\calH$, and a unit vector, $\varphi$, and the one parameter family
\begin{equation}\label{1.1}
  A_\alpha = A+\alpha\jap{\varphi,\cdot}\varphi
\end{equation}
If $F_\alpha(z) = \jap{\varphi,\tfrac{1}{A_\alpha-z}\varphi}$, then an easy calculation proves the \emph{Aronszajn-Krein formula}
\begin{equation}\label{1.2}
  F_\alpha(z) = \frac{F_0(z)}{1+\alpha F_0(z)}
\end{equation}
There is a wonderful set of spectral consequences of this formula that one can read about in Simon~\cite[Section 5.8]{OT} or Simon~\cite[Chaps 11-12]{TI}.  Included is the fact that the a.c. parts of the spectra of $A_\alpha$ (as $\alpha$ varies) are mutually absolutely continuous while the singular parts are mutually singular (when $\varphi$ is cyclic for $A$).

\bigskip

The Weyl-Titchmarsh $m$-function (defined below) is the basic object for the inverse theory of $1D$ (continuum) Schr\"{o}dinger operators (acting as an unbounded operator on $L^2(0,\infty)$ or on $L^2(\bbR)$)
\begin{equation}\label{1.3}
  Hu(x) = -u''(x)+V(x)u(x)
\end{equation}
In \eqref{1.3}, the potential, $V(x)$, can be taken as an $L^1_{\text{loc}}$ function but in this introduction for simplicity we'll suppose $V$ is continuous.  In the half line case, one needs to place a boundary condition at $x=0$ to define a self-adjoint operator (we discuss the situation at $\infty$ below).  We'll mainly fix Neumann conditions, i.e. $u'(0)=0$ in this introduction although we'll discuss more general boundary conditions in Section \ref{s2}.

One considers solutions of $Hu=zu$ for $z\in\bbC$ -- not as an operator statement but as a differential equation (i.e. we do not apriori ask that $u$ be in $L^2$).  For any $z\in\bbC_+$, the open upper half plane, there is a non-zero solution in $L^2$ at $+\infty$.  We say that $V$ is \emph{limit-point} (resp. \emph{limit circle}) at $\infty$ if for one, equivalently all, $z\in\bbC_+$, this $L^2$ solution is unique up to constant (resp. all of the $2D$ set of solutions are $L^2$).  It is a fundamental fact (see, for example, Simon~\cite[Theorem 7.4.12]{OT}) that in the $\bbR$ case, one has that $H$ defined initially on $C_0^\infty(\bbR)$ is essentially self-adjoint if and only $V$ is limit point at $\pm \infty$ and, in the $(0,\infty)$ case when $V$ is limit point at $\infty$, the operator on $C_0^\infty(0,\infty)$ has a one parameter family of self-adjoint extensions specified by a boundary condition of the form $au'(0)+bu(0)=0$ for fixed $a, b$ not both zero (this is in the case $V$ is continuous at $0$; the situation in the general $L^1_{\text{loc}}$ case is more involved).  Roughly speaking if $V(x) > -cx^\alpha$ with $\alpha<2$, and, in particular, if $V$ is bounded below, then $V$ is limit point at $\infty$ (see Simon~\cite[Theorems 7.4.19, 7.4.21]{OT}).

When $V$ is limit point at infinity, the \emph{Weyl-Titchmarsh $m$-function} is defined on $\bbC_+$ by (where $u_+$ is the solution which is $L^2$ at $\infty$)
\begin{equation}\label{1.4}
  m(z) = \frac{u_+'(0,z)}{u_+(0,z)}
\end{equation}
This is the m-function for Dirichlet BC.  The one for Neumann BC is
\begin{equation}\label{1.4A}
  \wti{m}(z) = -\frac{1}{m(z)}
\end{equation}
More generally, one sometimes defines
\begin{equation}\label{1.5}
  m(z,x) = \frac{u_+'(x,z)}{u_+(x,z)}
\end{equation}
which obeys the well-known \emph{Riccati equation} (with $m'=\tfrac{\partial m}{\partial x}$)
\begin{equation}\label{1.6}
  m'(z,x) = V(x) -z - m(z,x)^2
\end{equation}
In the $\bbR$ case, this is often written $m_+$.  There is then an $m_-$ in terms of the solution which is $L^2$ at $-\infty$.

A major result in inverse spectral theory (which Fritz and I often referred to) is

\begin{theorem} [Borg-Marchenko Uniqueness Theorem~\cite{Borg, MarchenkoPaper, MarchenkoPaper2}] \label{T1.1} Given $V_1, V_2$ on $[0,\infty)$, with m-functions $m_1, m_2$, if $m_1=m_2$ on all of $\bbC_+$, then $V_1=V_2$.
\end{theorem}

\begin{remark} At about the same time, Tikhonov~\cite{Tik} proved a special case of this result.
\end{remark}

While I called this a result in inverse spectral theory, it isn't clear from what was said so far why $m$ encapsulates spectral information.  In fact, when $H$ is defined by \eqref{1.3} with Neumann BC at $0$, the operator $(H-z)^{-1}$ for $z\in\bbC_+$ has a Green's function, i.e. $G(x,y;z)$ obeying
\begin{equation}\label{1.7}
  [(H-z)^{-1}u](x) = \int G(x,y;z)u(y)\,dy
\end{equation}
with (in case $V$ is continuous) $G$ jointly continuous on $[0,\infty)\times [0,\infty) \times\bbC_+$ and there is a measure, $d\rho$ on $\bbR$ with $\int (1+|\lambda|)^{-1}\,d\rho(\lambda)<\infty$ so that $H$ is unitarily equivalent to multiplication by $x$ on $L^2(\bbR,d\rho)$ and (for a suitable constant, $c$)
\begin{equation}\label{1.8}
  \wti{m}(z) = G(0,0;z) = \int \frac{1}{\lambda-z}\,d\rho(\lambda)
\end{equation}
One can go from the \emph{spectral measure}, $\rho$, to $m$ and by the theory of Stieltjes transforms (see Simon~\cite[Section 5.8]{OT}) from $m$ to $\rho$.  Thus, Theorem \ref{T1.1} is a result about the spectral measure determining the potential.  It is also a result in inverse scattering theory since when $V$ has sufficient decay, one can write the spectral measure in terms of the reflection coefficient (which determines the a.c. part of $\rho$) and the locations and weights of the bound states (which determine the pure point part of $\rho$).  This fact and a constructive way of recovering $V$ in this case is the content of the fundamental works of Gel'fand-Levitan~\cite{GelLev} and Marchenko~\cite{MarchenkoBk}.

\begin{remark} Technical point.  The condition $\int (1+|\lambda|)^{-1}\,d\rho(\lambda)<\infty$ needed for the representation \eqref{1.8} holds for Neumann and any BC except Dirichlet.  For Dirichlet one needs to replace $-1$ by $-2$ in the integral condition and the integral formula for $m$ takes the form \eqref{2.4}.

\end{remark}

Further insight can be obtained by considering the discrete case; indeed, Fritz and I wrote a partially pedagogical paper~\cite{Fritz24} emphasizing that the theory of orthogonal polynomials on the real line could be rephrased and understood in terms of an analogy to the inverse theory of Schr\"{o}dinger operators (that paper also has a number of new results obtained by using the analogy).  It should be emphasized that Harry Hochstadt~\cite{Hoch1, FuHoch, Hoch2, Hoch3} was a pioneer in exploiting this analogy.

One starts with a semi-infinite Jacobi matrix (we also study the two sided case on $\ell^2(\bbZ)$) with $b_1,b_2,\dots$ on diagonal and $a_1,a_2,\dots$ immediately off-diagonal, i.e.
\begin{equation} \lb{1.9}
  J(\{a_n,b_n\}_{n=1}^\infty)=\begin{pmatrix}
            b_1 & a_1 & 0  \\
            a_1 & b_2 & a_2 & \ddots \\
            0 & a_2 & b_3 & \ddots \\
            {} & \ddots & \ddots & \ddots & \ddots
    \end{pmatrix}
\end{equation}
as an operator on $\ell^2(\bbZ_+)$.  The parameters $\{a_n,b_n\}_{n=1}^\infty$ with $(a_n > 0; b_n\in\bbR)$ are called the \emph{Jacobi parameters}.  In the discussion that follows in this section, we will only consider the case where $\sup_n(|a_n|+|b_n|)<\infty$ so $J$ defines a bounded self-adjoint operator.

There is a unique \emph{spectral measure}, $d\mu$ on $\bbR$, so that
\begin{equation}\label{1.10}
  \jap{\delta_1,f(J)\delta_1} = \int f(\lambda)\,d\mu(\lambda)
\end{equation}
initially for polynomials, $f$, and then for all Borel functions (this is the spectral theorem for bounded operators; see Simon~\cite[Chapter 5]{OT}).  In particular, one has that
\begin{equation}\label{1.11}
  m(z) \equiv \bigjap{\delta_1,\frac{1}{J-z}\delta_1} = \int \frac{1}{\lambda-z}\,d\mu(\lambda)
\end{equation}
is the analog of the Weyl-Titchmarsh m-function of the continuum situation.  When we want to emphasize the Jacobi parameters involved we will write $m(z;\{a_n,b_n\}_{n=1}^\infty)$.

There is an analog of $m(z,x)$ denoted $m_n(z)$ obtained by stripping off the first $n$ rows and columns of the matrix, i.e. (we then sometimes write $m=m_0$)
\begin{equation}\label{1.12}
  m_n(z;\{a_k,b_k\}_{k=1}^\infty) \equiv m(z;\{a_k,b_k\}_{k=n+1}^\infty)
\end{equation}
The analog of the Riccati equation is
\begin{equation}\label{1.13}
  m_n(z) = \frac{1}{z-b_{n+1}+a_{n+1}^2m_{n+1}(z)}
\end{equation}

This non-trivial equation has several illuminating proofs (it suffices by induction to prove it when $n=0$).  One proof uses the fact that $m_1$ is a matrix element of a resolvent of $(1-P)J(1-P)$, where $J$ is the original Jacobi matrix, $P$ is the projection onto $\delta_1$ and we are interested in $P(J-z)^{-1}P$.  There is a general formula (called Banachiewicz' formula after \cite{Ban}) of relating those objects; this formula is associated in the mathematics literature with the method of Schur complements and in the physics literature with the Feshbach projection method (after \cite{Fesh, Schur}).  It is especially useful in studying Jacobi matrices on trees but can be specialized to this $1D$ case; one place to read the details is Avni et. al~\cite[(6.3)]{ABS}.

A second proof considers the analog (for $n=0$) of this equation where the semi-infinite Jacobi matrix is approximated by the $N\times N$ matrix, $J_N$, obtained by keeping the first $N$ rows and columns of \eqref{1.9}.  The $m$ function is then written as a ratio of determinants.  By expanding $\det(J_N-z)$ in minors in the first row, one gets a relation among this determinant and the two determinants obtained by stripping over the top row and leftmost column or two rows and two columns which leads to \eqref{1.13} in the limit as $N\to\infty$.  A version of this argument in a more general form appeared in Gesztesy-Simon~\cite{Fritz24}; see also Simon~\cite[Proposition 5.1]{SimonRatio}.

A third proof uses analogs of Weyl solutions using the finite difference equation (for $n\ge 1$; we pick some value of $a_0>0$)
\begin{equation}\label{1.14}
  a_nu_{n+1}+(b_n-z)u_n+a_{n-1}u_{n-1}=0
\end{equation}
Using $\sup_n |a_n|<\infty$, a Wronskian argument proves that this equation has at most one solution (up to a multiplicative constant) which is $\ell^2$ at infinity.  If $z\in\bbC_+$, one sees that
\begin{equation}\label{1.15}
  u_n = \left\{
          \begin{array}{ll}
            [(J-z)\delta_1]_n, & \hbox{ if } n\ge 1 \\
            -1/a_0, & \hbox{ if } n=0
          \end{array}
        \right.
\end{equation}
solves this equation, so for any $\ell^2$ solution, u, we see that
\begin{equation}\label{1.16}
  m(z) = -\frac{u_1}{a_0u_0}
\end{equation}
It follows that
\begin{equation}\label{1.17}
  m_1(z) = -\frac{u_2}{a_1u_1}
\end{equation}
so that \eqref{1.14} implies \eqref{1.13}. See Simon~\cite[Section 3.2]{SiSz} for details of this argument.

Before leaving the discussion of the this discrete case, we need to remark that the solution of \eqref{1.14} with initial conditions $u_0=0, u_1=1$ are $u_n=p_{n-1}(z)$ polynomials of degree $n$ in $z$ (by a simple inductive argument).  They are in fact the orthonormal polynomials for the spectral measure, $d\mu$, as follows from noting that $p_n(J)\delta_1=\delta_{n+1}$.  This immediately gives one way of recovering the Jacobi parameters from the spectral measure - form the orthonormal polynomials and look at their recursion coefficients.

There is, of course, a second method.  Iterating \eqref{1.13} gives a continued fraction expansion
\begin{equation}\lb{1.18}
m(z) = \cfrac{1}{b_1 -z -
               \cfrac{a_1^2}{b_2 - z -
                    \cfrac{a_2^2}{b_3 - z - \cdots}}}
\end{equation}
so one can form $m(z)$ from $d\mu$ via \eqref{1.11} and recover the Jacobi parameters from \eqref{1.18}.  Both approaches are associated with Jacobi although Jacobi only had the continued fractions for finite Jacobi matrices - the infinite version is associated with Markov and Stieltjes.  In particular, Stieltjes proved convergence of the continued fraction expansion in the upper half plane even for certain unbounded positive Jacobi operators.

\bigskip

The remainder of this paper will discuss four groups of results from my joint work with Fritz.  In section \ref{s2}, we focus on rank one perturbations, \eqref{1.1}, looking at a more general setting and discussing the limit of \eqref{1.1} as $\alpha\to\infty$.

Section \ref{s3} links the theory of rank one perturbations to inverse spectral theory.  In \cite{KreinSS}, Krein introduced a spectral shift function for operator differences that were trace class by building it up from a shift he could define for rank one perturbations.  In \cite{Fritz18}, we realized that the spectral shift, $\xi(x,\lambda)$ for the rank one perturbation that inserts a Dirichlet boundary condition at $x$ could be used to recover the potential via a trace formula that generalized many earlier special cases.

In Section \ref{s4}, we discuss the continuum analog of the Stieltjes continued fraction  approach to the inverse problem for Jacobi matrices. Finally, in Section \ref{s5}, we discuss our joint work that was motivated by a lovely result of Hochstadt-Lieberman~\cite{HochLieb} that for Schrodinger operators on $L^2((0,1),dx)$, the knowledge of the eigenvalues for the Neumann BC operator and of $V$ on $(0,\tfrac{1}{2})$ determines $V$ on all of $(0,1)$.  The recurrent theme in this work is how much information one needs about a Herglotz function in the upper half plane to determine it completely.

We close this section that focuses on inverse problems by remarking that we really know rather little about this subject.  This is illustrated by a problem that is open that Fritz and I have emphasized (for example in \cite{Fritz30} where we proved an analogous result for the half line): is the set of $C^\infty$ function on $\bbR$ whose spectrum is the same as the harmonic oscillator a connected set in some reasonable topology?

I want to thank Fritz Gesztesy, not only for the joy of collaboration but for useful feedback on the first draft of this article.  

%%%%%%%%%%%%%%%%%%%%%%%%%%%%%%%%%%%%%%%%%%%%%%%%%%%%%%%%%%%%%%
\section{Rank One Spectral Theory} \lb{s2}
%%%%%%%%%%%%%%%%%%%%%%%%%%%%%%%%%%%%%%%%%%%%%%%%%%%%%%%%%%%%%%

The m-function that we defined in Section \ref{s1} for the $1D$ continuum Schr\"{o}dinger operator, \eqref{1.3}, is appropriate for $u(0)=0$ Dirichlet BC.  More generally (e.g. Levitan-Sargsyan~\cite{LevSar}), one wants to consider for fixed $\theta\in [0,\pi)$, the boundary condition
\begin{equation}\label{2.1}
  u(0)\cos(\theta)+u'(0)\sin(\theta) = 0
\end{equation}
$\theta=0$ is Dirichlet BC which requires some special consideration while $\theta=\tfrac{\pi}{2}$ is Neumann BC.  For $\theta\ne 0$, one defines $m_\theta$ in terms of the Green's function for the operator, $H_\theta$, with BC \eqref{2.1} by
\begin{equation}\label{2.2}
  G_\theta(0,0;z) = \sin^2(\theta)[-\cot(\theta)+m_\theta(z)]
\end{equation}
There is a spectral measure, $d\rho_\theta$, so that $H_\theta$ is multiplication by $x$ in $L^2(\bbR,d\rho_\theta)$ and
\begin{equation}\label{2.3}
  m_\theta(z) = \cot(\theta) + \int_\bbR \frac{d\rho_\theta(\lambda)}{\lambda-z}
\end{equation}

For $\theta=0$, the Green's function vanishes at $x=0$ and one has to define $m$ in terms of derivatives of $G$ and the spectral representation takes the form
\begin{equation}\label{2.4}
 m(z) = c + \int \left[\frac{1}{\lambda-z}-\frac{\lambda}{1+\lambda^2}\right]\,d\rho(\lambda)
\end{equation}

One can relate $m_\theta$ to the Dirichlet  $m$ of \eqref{1.4} by (notice that \eqref{1.4A} is just this equation for $\theta=\tfrac{\pi}{2}$)
\begin{equation}\label{2.5}
  m_\theta(z) = \frac{\cos(\theta)m(z)-\sin(\theta)}{\sin(\theta)m(z)+\cos(\theta)}
\end{equation}
This has a similarity to \eqref{1.2} and indeed, spectral consequences of \eqref{1.2} can be mimicked to say the same things for variation of BC. In fact, Aronszajn~\cite{Ar3} proved what we listed as consequences of the rank one theory for variation of BC and Donoghue~\cite{Donog} used the similarity of formulae to extend them to the rank one setting.

Notice one big difference between these similar settings: the parameter $\theta$ is naturally a circle because $0$ and $\pi$ are the same while the parameter space for rank one perturbations is $\bbR$.

One of the points of the paper of Gesztesy-Simon\cite{Fritz15} which is the subject of this section was to popularize and probably to initiate a framework for the theory of rank one perturbations where change of BC was not merely an analog - rather it became a special case (at least in the situation where all the operators are bounded from below).  One needs to allow $A$ to be unbounded but we only considered the case where $A$ was bounded from below.  In fact we only considered $A\ge 0$.  The perturbation was rank one in a certain sense but not required to be bounded.

We recall the basic facts about quadratic forms associated to positive self-adjoint operators, $A$, on a Hilbert space, $\calH$, as discussed, for example, in Simon~\cite[Section 7.5]{OT}.  One forms the form domain, $Q(A)$, those vectors, $\varphi\in\calH$, for whose spectral measures, $d\mu_\varphi$, obey
\begin{equation}\label{2.6}
  \int \lambda \,d\mu_\varphi(\lambda) < \infty
\end{equation}
in which case the integral is formally denoted as $\jap{\varphi,A\varphi}$ (which is what it is if it happens that $\varphi\in D(A)$).

One puts a norm, $\Norm{\cdot}_{+1}$, on $Q(A)$ by
\begin{equation}\label{2.7}
  \Norm{\varphi}_{+1}^2 = \jap{\varphi,A\varphi}+\Norm{\varphi}^2
\end{equation}
Then $Q(A)$ is complete in this norm and it comes from an inner product so we call this Hilbert space $\calH_{+1}$.  We let $\calH_{-1}$ be its dual space (where we suppress that a Hilbert space can be identified with its dual).  If $\psi\in\calH$, we can define a linear functional on $\calH_{+1}$, by $\varphi\in\calH_{+1}\mapsto\jap{\psi,\varphi}$ (with the $\calH$, not $\calH_{+1}$ inner product). For this reason, we denote the duality of $\calH_{-1}$ on $\calH_{+1}$ by $\varphi(\eta)\equiv\jap{\varphi,\eta}$ using $\jap{\eta,\varphi}$ for the complex conjugate.  In this way, we get dense inclusions
\begin{equation}\label{2.8}
  \calH_{+1}\subset\calH=\calH^*\subset\calH_{-1}
\end{equation}
This scale of subspaces is a standard construction associated with positive self-adjoint operators.

Now suppose $\varphi\in\calH_{-1}$ (which we assume is cyclic for $A$ in the sense that the linear span of $\{(A-z)^{-1}\varphi\,\mid\,z\in\bbC_+\}$ is dense in $\calH$).  Then $\eta\in\calH_{+1}\mapsto \jap{\eta,\varphi}\jap{\varphi,\eta}$ defines a quadratic form with form domain $Q(A)$.  We denote it as $\jap{\varphi,\cdot}\varphi$.  It is easy to see that in the language of \cite[Theorem 7.5.7]{OT} this form is a form bounded form perturbation of $A$ of relative bound $0$, so for all $\alpha\in\bbR$, \eqref{1.1} (where now $\varphi$ might no longer be a unit vector) defines a one parameter family of closed semibounded quadratic forms and so a one parameter family of self-adjoint operators.  If $A$ is bounded and $\varphi\in\calH$, it is exactly that class discussed after \eqref{1.1}.  If $V$ is continuous and bounded from below, one can take $A$ to be the Neumann BC operator, $H_{\theta=\tfrac{\pi}{2}}$.  One can show (via a Sobolev estimate) that $Q(A)$ consists of continuous functions and $\eta\mapsto \eta(0)$ lies in $\calH_{-1}$, a vector denoted $\delta(0)$.  It is easy to see that for $\theta\ne 0,\pi$, we have that
\begin{equation}\label{2.9}
  H_\theta = A-\cot(\theta)\jap{\delta(0),\cdot}\delta(0)
\end{equation}
so BC variation is a special case of rank one perturbation.

While the formalism and BC as a special case are important parts of \cite{Fritz15}, as the title of the paper, \emph{Rank one perturbations at infinite coupling}, suggests, the main point is to show  that the parameter space is not $\bbR$ but $\bbR\cup\{\infty\}$, the one point compactification of $\bbR$ thereby resolving the circle/real line puzzle noted above.  It is a general fact (see Simon~\cite[Theorem 7.5.14]{OT}) that if $0\le A\le B$ as quadratic forms (in the sense that $Q(B)\subset Q(A)$ and, for $\varphi\in Q(B)$, one has that $\jap{\varphi,A\varphi}\le\jap{\varphi,B\varphi}$), then for all $x\ge 0$, one has that $(B+x)^{-1}\le (A+x)^{-1}$.

For $A_\alpha$ given by \eqref{1.1} in this general setting, one thus sees that as $\alpha\to +\infty$, $(A_\alpha+x)^{-1}$ has a strong limit (since its expectations are monotone).  A general result known as the monotone convergence theorem for forms (see Kato~\cite[Section VIII.2]{KatoMC} and Simon~\cite{SimonMC}) identifies the limiting operator.  Let
\begin{equation}\label{2.10}
  Q(A_\infty) = \{\psi\in\cap_\alpha Q(A_\alpha)\,\mid\, \sup_\alpha \jap{\psi,A_\alpha\psi} <\infty\}
\end{equation}
For $\psi\in Q(A_\infty)$, one defines
\begin{equation}\label{2.11}
  \jap{\psi,A_\infty\psi} = \sup_\alpha \jap{\psi,A_\alpha\psi}
\end{equation}
This is a closed quadratic form but it may not be densely defined.  If $\calH^{(\infty)}$ is the closure of $Q(A_\infty)$ in $\calH$, the general theory of closed quadratic forms (see Simon~\cite[Theorem 7.5.5]{OT}) says that there is a self-adjoint operator on $\calH^{(\infty)}$ so that if $P:\calH\to\calH^{(\infty)}$ is the orthogonal projection onto $\calH^{(\infty)}$, then one has that for all $\eta\in\calH$ and $z\in\bbC_+$
\begin{equation}\label{2.12}
  \lim_{\alpha\to\infty}(A_\alpha-z)^{-1}\eta = (A_\infty-z)^{-1}P\eta
\end{equation}

If one specializes this monotone convergence theorem to the rank one situation of \eqref{1.1} with $\varphi\in\calH_{-1}$ one sees that $\calH^{(\infty)}=\calH$ if $\varphi\in\calH_{-1}\setminus\calH$ while if $\varphi\in\calH$, then $\calH_\infty$ is the orthogonal complement of $\{\varphi\}$. In both cases, $Q(A_\infty)$ is strictly smaller than $\calH_{+1}$, indeed, it is $\{\eta\in\calH_{+1}\,\mid\,\jap{\varphi,\eta}=0\}$.

Beyond exploiting the monotone convergence theorem, Gesztesy-Simon~\cite{Fritz15} use the explicit formula
\begin{equation}\label{2.13}
  (A_\alpha-z)^{-1} = (A-z)^{-1} - \frac{\alpha}{1+\alpha F_0(z)}\jap{(A-\bar{z})^{-1}\varphi,\cdot}(A-z)^{-1}\varphi
\end{equation}
which has the advantage of showing that when $z\in\bbC_+$, the strong limit is the same as $\alpha\to -\infty$ as when $\alpha\to +\infty$.

It follows from \eqref{1.2} that $\lim_{\alpha\to\infty} F_\alpha(z)=0$ so that if $d\mu_\alpha$ is the spectral measure for $A_\alpha, \varphi$, then for all continuous functions, $f:\bbR\to\bbR$, of compact support, one has that
\begin{equation}\label{2.14}
  \lim_{\alpha\to\infty} \int f(\lambda)\,d\mu_\alpha(\lambda) = 0
\end{equation}

However, Gesztesy-Simon~\cite{Fritz15} showed that if
\begin{equation}\label{2.15}
  d\rho_\alpha(\lambda) = (1+|\alpha|^2)d\mu_\alpha(\lambda)
\end{equation}
then there is a non-zero measure $d\rho_\infty$ so that for continuous functions, $f:\bbR\to\bbR$, of compact support one has that
\begin{equation}\label{2.16}
  \lim_{\alpha\to\infty} \int f(\lambda)\,d\rho_\alpha(\lambda) = \int f(\lambda)\,d\rho_\infty(\lambda)
\end{equation}
Moreover, there is a vector, $\eta\in\calH_{-2}(A_\infty)$ which is cyclic for $A_\infty$ so that $d\rho_\infty$ is the spectral measure for $A,\eta$.   In addition, they prove that for such $f$
\begin{equation}\label{2.17}
  \int f(\lambda)\,d\rho_\infty(\lambda) = \lim_{\varepsilon\downarrow 0} \int f(\lambda)\Ima\left[\frac{-1}{F_0(\lambda+i\varepsilon)}\right]\,d\lambda
\end{equation}

If one looks at what this says for the situation of boundary conditions given in \eqref{2.9}, then $A_\infty$ is the Dirichlet BC operator and \eqref{2.17} is just \eqref{1.4A}.  Thus the $A_\infty$ construction of Gesztesy-Simon~\cite{Fritz15} is just the abstraction of the construction of Dirichlet BC.

Among the papers that followed up on this work, we should mention that there is an interesting collection that generalizes rank one perturbations to allow $\varphi\in\calH_{-2}$, see Kiselev-Simon~\cite{KisSi} and Albeverio et. al.~\cite{AKK}. We also mention two papers of Hassi et al.~\cite{Hassi1, Hassi2} one of which allows $A$'s which are unbounded above and below and the other explores what happens when $\varphi\in\calH_{+k}$ with $k\ge 1$.

%%%%%%%%%%%%%%%%%%%%%%%%%%%%%%%%%%%%%%%%%%%%%%%%%%%%%%%%%%%%%%
\section{The Xi Function} \lb{s3}
%%%%%%%%%%%%%%%%%%%%%%%%%%%%%%%%%%%%%%%%%%%%%%%%%%%%%%%%%%%%%%

In this section, we will discuss a series of papers by Gesztesy and Simon, some with Holden and/or Zhao~\cite{Fritz13, Fritz18, Fritz14, Fritz16, Fritz17} on trace formulae of which \cite{Fritz13} is an announcement and \cite{Fritz18} the central one.  While there are some results for higher dimensional Schr\"{o}dinger operators~\cite{Fritz17}, for higher order (i.e. for derivatives of $V$)~\cite{Fritz16} and for Jacobi matrices~\cite{Fritz13}, I will restrict the discussion here to the lowest order trace formulae for continuum Schr\"{o}dinger operators on the line.

Like the m-function, the basic function will be connected to diagonal Green's function (integral kernel, $G(x,y;z)$) of \eqref{1.3} but on the whole line.  The \emph{xi function}, $\xi(x,\lambda)$, is given by
\begin{equation}\label{3.1}
  \xi(x,\lambda) = \frac{1}{\pi} \arg\left(\lim_{\varepsilon\downarrow 0}G(x,x;\lambda+i\varepsilon)\right)
\end{equation}
where the general theory of Herglotz functions (see Simon~\cite[Section 5.9]{HA}) implies that for each fixed $x$, the limit exists for Lebesgue a.e. $\lambda$ and that
\begin{equation}\label{3.2}
  0\le \xi(x,\lambda)\le 1
\end{equation}

The basic trace formula says that for any $V$ bounded from below and any $E_0\le\inf\spec(H)$, one has that
\begin{equation}\label{3.3}
  V(x) = \lim_{\alpha\downarrow 0}\left[E_0+\int_{E_0}^{\infty}e^{-\alpha\lambda}[1-2\xi(x,\lambda)]\,d\lambda\right]
\end{equation}
\begin{remarks} 1. We intend this for continuous $V$ at all $x$ but one can prove it (\cite{Fritz16}) for $V\in L^1_{\text{loc}}$ for a.e. $x$.  It was proven in \cite{Fritz18} for continuous $V$ obeying $V(x) \le C_1\exp(C_2x^2)$ for some $C_1, C_2>0$ and for all $x$ and in \cite{Fritz16} for all $V$ bounded from below without restriction on growth.

2. Below $\spec(H)$, $\xi\equiv 0$ so the right side of \eqref{3.3} is unchanged if $E_0$ is changed so long as $E_0\le\inf\spec(H)$.

3. The limit in \eqref{3.3} can be viewed as an abelian summation. If $H$ has discrete spectrum (e.g. if $V\to\infty$ as $x\to\pm\infty$), then $|1-2\xi(x,\lambda)|=1$ for all $x$ and a.e. $\lambda$ so the integral in \eqref{3.3} is not absolutely convergent when $\alpha=0$ and the summation method is needed.

4.  Of course, if $\int_{E_0}^{\infty}|1-2\xi(x,\lambda)|\,d\lambda<\infty$, by the dominated convergence theorem, we can take $\alpha=0$ in the integral and drop the $\lim$.

5.  The term \emph{trace formula} is overused in mathematics - however, not only is this a generalization of what others have called a trace formula but, as we will see, it does come from the calculation of a suitable trace!
\end{remarks}

\bigskip

Given our length limitations we cannot give complete proofs but will say something about the history, the interpretation of $\xi$ as a Krein spectral shift, make some remarks on the proof and give several examples. The earliest trace formulae for Schr\"{o}dinger operators were found on a finite interval in 1953 by Gel'fand and Levitan~\cite{GelLevTrace} with later contributions by Dikii~\cite{DikTrace}, Gel'fand~\cite{GelTrace}, Halberg-Kramer~\cite{HKTrace}, and Gilbert-Kramer~\cite{GKTrace}.

After that there was a lot of work on periodic Schrodinger operators, so suppose that
\begin{equation}\label{3.4}
  V(x+1) = V(x)
\end{equation}
The beautiful theory of the basic structure of the spectrum is described in Magnus-Winkler~\cite{MW}.  There are a sequence of real numbers (the eigenvalues of the operators on $L^2([0,2], dx)$ with periodic BC written in increasing order counting multiplicity): $E_0<E_1\le E_2<E_3\le E_4\dots E_{2j-1}\le E_{2j}< E_{2j+1}$.  The spectrum of full line operators is
\begin{equation}\label{3.5}
  \spec(H) = [E_0, E_1]\cup\dots\cup[E_{2j},E_{2j+1}]\cup\dots
\end{equation}
so there are gaps in the spectrum (where, it is known, the spectral measures are purely absolutely continuous).  The gaps tend to shrink with the shrinkage related to smoothness of $V$.  In particular, if $V$ is $C^1$, then
\begin{equation}\label{3.6}
  \sum_{j=1}^{\infty} |E_{2j} - E_{2j-1}| < \infty
\end{equation}

Consider for each $y\in\bbR$, the operator, $H_y$, on $L^2([y,y+1],dx)$ with Dirichlet BC, $u(y)=u(y+1)=0$.  It has an infinite set of eigenvalues, $\mu_1(y)<\mu_2(y)<\dots$ which are known to obey (for each $y$)
\begin{equation}\label{3.7}
  E_{2j-1} \le \mu_j(y) \le E_j
\end{equation}
What the earlier trace formula says is that when $V$ is $C^1$, we have that
\begin{equation}\label{3.8}
  V(y) = E_0+\sum_{j=1}^{\infty} [E_{2j-1}+E_{2j}-2\mu_j(y)]
\end{equation}
By \eqref{3.7} we have that (because the left side is twice the distance from $\mu_j(y)$ to the midpoint of the interval)
\begin{equation}\label{3.9}
  |E_{2j-1}+E_{2j}-2\mu_j(y)| \le |E_{2j}-E_{2j-1}|
\end{equation}
so the sum in \eqref{3.8} is absolutely convergent by \eqref{3.6}.

The first formula related to this for some periodic $V$'s was found in 1965 by Hochstadt~\cite{HochTrace}, namely, he proved that in the finite gap case that
\begin{equation}\label{3.10}
  V(y) - V(0) = 2\sum_j  [\mu_j(y)-\mu_j(0)]
\end{equation}
Dobrovin~\cite{DobTrace} then proved \eqref{3.8} for the finite gap case and McKean-van Moerbeke~\cite{McvMTrace} and Flaschka~\cite{FlTrace} proved \eqref{3.8} for sufficiently smooth periodic $V$'s.

Before leaving the history, we want to mention one other prior result that is a trace formula.  Namely, in 1979, Deift-Trubowitz~\cite{DTTrace} for a whole line problem with no bound states and sufficiently rapid decay on $V$, one has that
\begin{equation}\label{3.11}
  V(x) = \frac{2i}{\pi} \int_{-\infty}^{\infty} k \log\left[1+R(k)\frac{f_+(x,k)}{f_-(x,k)}\right]\,dk
\end{equation}
where $R(k)$ is the reflection coefficient and $f_\pm(x,k)$ the Jost functions for energy $E=k^2$.  Remarkably \eqref{3.8} and \eqref{3.11} are special cases of the same formula, namely \eqref{3.3}!

\bigskip

The $\xi$-function is actually a Krein spectral shift so we need to quickly recall what that it is, especially in the rank one situation, \eqref{1.1}.  The general notion (Simon~\cite[Appendix to Section 11.4]{TI}) is that given a self-adjoint operator, $A$ bounded from below and $C$  trace class (see Simon~\cite{TI} for the notion of trace class) with $B=A+C$, there is a unique $L^1$ function, $\xi_{B,A}(x)$, the \emph{Krein spectral shift}~\cite{KreinSS}, on $\bbR$ determined by the conditions that $\xi_{B,A}$ vanishes near $-\infty$, and so that for all $C^2$ functions with $f, f', f''$ bounded by a multiple of $(1+|x|^2)^{-1}$ one has that $f(B)-f(A)$ is trace class and
\begin{equation}\label{3.12}
  \tr(f(B)-f(A)) = \int f'(x)\xi_{B,A}(x)\,dx
\end{equation}

One way of approaching the construction of this function (see Simon~\cite[Sections 5.8-5.9]{OT}) is to do it for rank one and use that to build up for general trace class perturbations proving along the way that
\begin{equation}\label{3.13}
  \int |\xi_{B,A}(x)|\,dx \le \Norm{B-A}_1
\end{equation}
For this rank one case with $A_\alpha$ given by \eqref{1.1}, one can write $\xi_\alpha \equiv \xi_{A_\alpha,A}$ explicitly (see Simon~\cite[Section 11.4]{TI})
\begin{equation}\label{3.14}
  \xi_\alpha(\lambda) = \left\{
                         \begin{array}{ll}
                           \tfrac{1}{\pi}\arg(1+\alpha F_0(\lambda+i0)), & \hbox{ if }|\alpha|<\infty  \\[3pt]
                           \tfrac{1}{\pi}\arg(F_0(\lambda+i0)), & \hbox{ if }\alpha=\infty
                         \end{array}
                       \right.
\end{equation}

\bigskip

Looking at the rank one formulation of adding a Dirichlet BC, \eqref{2.9} and \eqref{3.14}, one sees that $\xi(x,\cdot)$, defined in \eqref{3.1}, is just the Krein spectral shift for going from $H$ to $H_{D;x}$ where the latter is the operator on $L^2((-\infty,x),dx)\oplus L^2((x,\infty),dx)$ with a Dirichlet BC at $x$ on both half line pieces. In particular, we have, by \eqref{3.12}, that
\begin{equation}\label{3.15}
  \tr(e^{-tH}-e^{-tH_{D;x}}) = t\int_{E_0}^{\infty} e^{-t\lambda}\xi(x,\lambda)\,d\lambda
\end{equation}
The other fact one needs to get to the general trace formula, \eqref{3.3}, is small time asymptotics when $x$ is a point of Lebesgue continuity for $V$
\begin{equation}\label{3.16}
  \tr(e^{-tH}-e^{-tH_{D;x}}) = \tfrac{1}{2}[1-tV(x)+\text{o}(t)]
\end{equation}
which is shown using path integrals in increasing degrees of generality in \cite{Fritz18, Fritz16, Fritz14}.  To get \eqref{3.3} from the last two formulae one subtracts the integral in \eqref{3.15} from
\begin{equation}\label{3.17}
  \tfrac{1}{2}\alpha\int_{E_0}^{\infty}e^{-\alpha\lambda}\,d\lambda = \tfrac{1}{2}[1-\alpha E_0+\text{0}(\alpha)]
\end{equation}

\bigskip

We close this section with a summary of examples:

(1) In the periodic case, it is known (e.g. Kotani~\cite{Kotani}), that $H$ is reflectionless, i.e. $G(x,x;\lambda+i0)$ is pure imaginary on the spectral bands, and that $G(x,x;\lambda)$ is real and monotone in each gap with a zero exactly at $\mu_j(x)$, thus for $\lambda\ge E_0$, we have that
\begin{equation}\label{3.18}
  \xi(x,\lambda) = \left\{
                     \begin{array}{ll}
                       \tfrac{1}{2}, & \hbox{ if } E_{2j}<\lambda<E_{2j+1}\\
                       1, & \hbox{ if } E_{2j+1}<\lambda<\mu_{j+1}(x)\\
                       0, & \hbox{ if } \mu_{j+1}(x)<\lambda<E_{2j+2}
                     \end{array}
                   \right.
\end{equation}
If $V$ is $C^1$ so \eqref{3.6} holds, then $\int_{E_0}^{\infty} |1-2\xi(x,\lambda)|\,d\lambda<\infty$ so we can take the limit inside the integral and we prove sum rule \eqref{3.8}.  Moreover, if \eqref{3.6} fails, we still have a result!

(2) Suppose $V(x)\to\infty$ as $x\to\pm\infty$ so $H$ has an infinite simple purely discrete spectrum with eigenvalues $E_0<E_1<\dots$. For each fixed $x$, $G(x,x;\lambda)$ is monotone in $\lambda$ in each interval between eigenvalues, running from $-\infty$ to $\infty$ with a zero at $\mu_j(x)$, the Dirichlet eigenvalue, which it is traditional to label starting with $\mu_1(x)$.  Thus $2\xi(x,\lambda)-1$ is $1$ on $[E_{j-1},\mu_j(x))$ and $-1$ on $(\mu_j(x),E_j]$, so we find the sum rule (new in \cite{Fritz13})
\begin{equation}\label{3.19}
  V(x) = E_0 + \lim_{\alpha\downarrow 0} \alpha^{-1}\sum_{j=1}^{\infty}[2e^{-\alpha\mu_j(x)}-e^{-\alpha E_j}-e^{-\alpha E_{j-1}}]
\end{equation}
An amusing special case is $V(x)=x^2-1$, so $\{E_j\}_{j=0}^\infty$ is $0,2,4,6,\dots$ and $\{\mu_j(0)\}_{j=0}^\infty$ is $2,2,6,6,10,10\dots$.  The integral $\int_{0}^{\infty}(1-2\xi(0,\lambda))\,d\lambda$ is $-2+2-2+2-2+\dots$ and \eqref{3.19} is the Abelian sum which is $-1=V(0)$.

(3) The paper of Gesztesy, Holden and Simon~\cite{Fritz14} discusses situations where is there is a.c. spectrum which is not reflectionless so regions where $0<\xi(x,\lambda)<1$ but so that $2\xi(x,\lambda)-1$ is not identically zero which means the a.c. spectrum contributes to the sum rule.  When $V$ decays to zero sufficiently fast at $\pm\infty$ and there are no bound states, an explicit calculation of $\xi$ in terms of the Jost solutions recovers the sum rule,\eqref{3.11}, of Deift-Trubowitz.  There are also calculations of the contributions of bound states, examples where $V$ approaches different constants (sufficiently quickly) at $\pm\infty$ and of a periodic potential with a short range perturbation (modelling impurities in crystals).

(4) Gesztesy-Simon~\cite[Theorem 5.5]{Fritz18} note the general fact (that follows from the fact that the ac part of the spectral measure associated to the $\calH_{-1}$ vector $\delta_x$ is $\Ima(G(x,x;\lambda+i0))\,d\lambda$) that the spectrum of the a.c. part of $H$ is, for each fixed $x$, the essential closure of the set of $\lambda$ where $0<\xi(\lambda,x)<1$.  They could use this to prove that if $V_n(x)\to V_\infty(x)$ uniformly on compact intervals with $\sup_{x,n}|V_n(x)|<\infty$, then for any Borel set $S\subset\bbR$, one has that (with $|\cdot|$ Lebesgue measure)
\begin{equation}\label{3.20}
  |S\cap\sigma_{ac}(H)| \ge \limsup_n |S\cap\sigma_{ac}(H_n)|
\end{equation}
This in turn has a consequence about which there is an interesting story concerning the \emph{almost Mathieu operator} (acting on $\ell^2(\bbZ)$).
\begin{equation}\label{3.21}
 (H_{\alpha,\lambda}(\theta)u)_n = u_{n+1}+u_{n-1}+\lambda\cos(\pi\alpha n+\theta)u_n
\end{equation}
where $\alpha,\lambda\in\bbR$ and $\theta\in [0,2\pi)$ (one normally only takes $\alpha\in [0,2)$ since it is invariant if $\alpha\to\alpha+2$).  We write $\theta$ as a parameter; one normally imagines looking at the whole family of operators for fixed $\alpha\in\bbR\setminus\bbQ,\lambda$ but varying $\theta$.  There are close connections as one varies $\theta$; for example the spectra are the same but the spectral types are only the same for (Lebesgue) a.e. $\theta$, not all $\theta$ in general.  From the work of Avron-Simon~\cite{AvSi}, it was known that for $\lambda>2$ and $\alpha$ an irrational very well approximated by rationals, that $H_{\alpha,\lambda}(\theta)$ has purely singular spectrum for all $\theta$ and this was very widely expected to hold true also for $|\lambda|<2$.  In May 1992, Yoram Last, a then graduate student of Yosi Avron, submitted a paper to me as an editor of CMP claiming that for \emph{all} $\alpha$ , the measure of the a.c. essential spectrum of $H_{\alpha,\lambda}(\theta)$ was at least $4-2|\lambda|$ and so there was such a.c. spectrum for \emph{all} $\alpha$.  This is contrary to the accepted expectation so I was sure it was wrong but the arguments were involved and I didn't see the error. Since I expected to be visiting Technion within a few weeks I conferred with Yosi to confirm that he and his student were OK to holding off on sending the paper to a referee until I could sit down with Yoram.  I told Yosi I was sure I'd find an error.  In fact using that in 1990 Avron, Simon  and van Mouche~\cite{AvSiVM} had proven
\begin{equation}\label{3.22}
  \alpha\in\bbQ\Rightarrow|\spec(H_{\alpha,\lambda}(\theta))|\ge 4-2|\lambda|
\end{equation}
Last had proven that the same inequality holds for all $\alpha$ and he convinced me that his argument was correct~\cite{LastMeas1, LastMeas2}.  Later that summer I returned to Caltech and worked with Fritz on the $\xi$-function and we realized \eqref{3.20}.  Moreover, this immediately implied that the result of  Avron, Simon  and van Mouche~\cite{AvSiVM} for rational $\alpha$ proved \eqref{3.22} for all real $\alpha$.  So, with help from Fritz, in a few weeks I went from a skeptic on the Last result to someone with a quick alternate proof!  There is, of course, much more known about this model now; see, for example Jitomirskaya~\cite{Jit1, Jit2}.

%%%%%%%%%%%%%%%%%%%%%%%%%%%%%%%%%%%%%%%%%%%%%%%%%%%%%%%%%%%%%%
\section{The A-Function and A Local Borg–Marchenko Uniqueness Theorem} \lb{s4}
%%%%%%%%%%%%%%%%%%%%%%%%%%%%%%%%%%%%%%%%%%%%%%%%%%%%%%%%%%%%%%

In this section, I want to discuss the A-function where Fritz and I wrote an important paper~\cite{Fritz28} which appeared in Ann. Math., a followup paper on an earlier paper that I wrote~\cite{SiInv}.  I will also briefly discuss a related paper that Fritz and I wrote~\cite{Fritz29}.  To set the stage, I will need to discuss the first paper~\cite{SiInv}.  Since the two parts total more than 80 pages, I will only hit the high points.

While \cite{SiInv} is not my most important or deepest paper, it is in many ways the one I am proudest of.  I tend to not be persistent.  If I work on a problem without success,  I usually not only drop the problem but don't return to it.  On this paper, I had a promising approach and for about ten years, I'd return to it periodically until I finally solved it.  The question concerns inverse spectral theory of continuum $1D$ Schrodinger operators.  As noted above for the discrete, i.e. Jacobi, case, there are two approaches for recovering the Jacobi parameters from the $m$-function: first, one can form orthogonal polynomials and get the Jacobi parameters from the difference equation that they obey and second, one can write down the continued fraction expansion, \eqref{1.18}, and read off the Jacobi parameters from that.

The continuum analog of the first approach is that of Gel'fand-Levitan~\cite{GelLev}, who followed the orthogonal polynomials route.  Their analog of monic polynomials were functions of the form ($\tfrac{\sin(kx)}{k}$ are the eigenfunctions for $H_0=-d^2/dx^2$ on the half line with $u(0)=0$)
\begin{equation}\label{4.1}
  s(x,k) = \frac{\sin(kx)}{k} + \int_{-x}^{x} K(x,y) \frac{\sin(ky)}{k} \,dy
\end{equation}
If $d\rho(\lambda)$ was the spectral measure for the Dirichlet BC operator with potential, $V$, given by \eqref{2.4}, they looked for $s(x,k)$ that obeyed \eqref{4.1} and
\begin{equation}\label{4.2}
  \int s(x,k)s(y,k) \,d\rho(k^2) = \delta(x-y)
\end{equation}
at least formally.  Once they constructed these analogs of the monic polynomials, they could recover $V$ from the fact that $s(x,k)$ solves $-u''+Vu=k^2 u$ (actually, they got it from intermediate objects in their theory).  But, as they remarked in their paper, they were motivated by following the recovery of Jacobi coefficients via orthogonal polynomials.

To be more explicit about the second approach to the Jacobi inverse problem, one begins with \eqref{1.13} for $n=0$ (with $m_0=m$).  From the leading behavior of $1/m(z)$ near $z=\infty$, we get $b_0$ and $a_0$ and then $m_1$, after which we can iterate.  I tried to find the continuum analog.  The initial step of reading $V(0)$ from $m$ is easy.  Atkinson~\cite{Atk} has proven that as $\kappa\to \infty$ (many of these hold in the complex $\lambda=-\kappa^2$ plane with a sector about the positive real axis removed), one has that
\begin{equation}\label{4.3}
  m(-\kappa^2) = -\kappa - \int_{0}^{a} V(\alpha)e^{-2\alpha\kappa}\,d\alpha + \text{o}(\kappa^{-1})
\end{equation}
for all small $a$ which implies (see also Danielyan-Levitan~\cite{DL}) that when $V$ is continuous that
\begin{equation}\label{4.4}
  m(-\kappa^2,x) = -\kappa - \frac{V(x)}{2\kappa} + \text{o}(\kappa^{-1})
\end{equation}
so that $m(\lambda,x)$ determines $V$.  But one wants to determine $V$ from $m(\lambda,x=0)$. Unlike the Jacobi case, where the iteration equation decouples from future Jacobi parameters, the differential equation for $m$ requires one to know $V$ not just at zero but at all small $x$!

In \cite{SiInv}, I managed the decoupling for two situations: first, where one considers an operator on $L([0,b],dx])$ for some $b<\infty$ with $u(0)=0$ and $u'(b)+hu(b)=0$ (for some fixed $h\in\bbR\cup\{\infty\}$) (where $h=\infty$ is short for $u(b)=0$ BC) boundary conditions and $\int_{0}^{b} |V(y)|\,dy <\infty$; second, where one takes $b=\infty$ and requires $\int_{0}^{a}|V(y)|\,dy < \infty$ for all $a<\infty$ and that $\sup_a \int_{a}^{a+1} \max(-V(y),0)\,dy <\infty$ (in which case it is known that no BC is needed at infinity). In this case, I was able to show that, for all $x$ and all $a<b$ and $\delta>0$, there was a representation
\begin{equation}\label{4.5}
  m(-\kappa^2,x) = -\kappa - \int_{0}^{a} A(\alpha,x) e^{-2\alpha\kappa}\,d\alpha + \text{O}(e^{-2\alpha\kappa+\delta})
\end{equation}
for a suitable locally $L^1$ function, $A$, of two variables.  $A(\cdot,x)$ is uniquely determined by $m(\cdot,x)$.  Moreover (and this is the version from~\cite{Fritz28}; \cite{SiInv} has somewhat weaker results), uniformly for $x$ in compacts of $[0,b)$, one has that
\begin{equation}\label{4.6}
  \lim_{\alpha\downarrow 0} |A(\alpha,x)-V(\alpha+x)|=0
\end{equation}

The final ingredient that makes everything work is the translation of the Riccati equation, \eqref{1.6} to an integro-differential equation
\begin{equation}\label{4.7}
  \frac{\partial A}{\partial x}(\alpha,x) = \frac{\partial A}{\partial\alpha}(\alpha,x) + \int_{0}^{\alpha} A(\alpha-\beta,x)A(\beta,x)\, d\beta
\end{equation}
which is proved to hold in distributional sense and in stronger senses when $V$ is smooth. Remarkably, $V$ has dropped out! Because of characteristics of the differential equation part, one can prove that $A(\alpha,x)$ in the region $\calR_K=\{(\alpha,x)\,\mid\,0<\alpha,x; 0<\alpha+x<K\}$ only depends on the values of $A(\alpha,x=0)$ for $0<\alpha<K$.  The solution of the inverse problem is thus clear - use \eqref{4.5} for $x=0$ to determine $A(\alpha,x=0)$ from $m(-\kappa^2)$, then \eqref{4.7} to determine $A(\alpha,x)$ and then \eqref{4.6} to determine $V$.  In particular, $A(\alpha,x=0)$ for $0<\alpha<K$ determines $V(x)$ for $0<x<K$.

A technical argument in the appendix of \cite{SiInv} proves that conversely, $V(x)$ for $0<x<K$ determines $A(\alpha,x=0)$ for $0<\alpha<K$.  In particular, this proves

\begin{theorem} [Local Borg-Marchenko Uniqueness Theorem]  \label{T4.1} Given two potentials, $V_1$ and $V_2$ with $m$-functions, $m_1$ and $m_2$, we have that $V_1(x)=V_2(x)$ for all $x\in [0,a]$ if and only if for all $\delta>0$, we have that as $\kappa\to\infty$
\begin{equation}\label{4.8}
  |m_1(-\kappa^2)-m_2(-\kappa^2)| = \text{O}(e^{-2a\kappa+\delta})
\end{equation}
\end{theorem}

\begin{remarks} 1. We have been sloppy about the conditions on $V_{1,2}$ because while \cite{SiInv} proves this for restricted $V$'s, \cite{Fritz28} proves it for all $L^1_{\text{loc}}$ potentials with any boundary conditions.

2. This has two directions whose proofs are different.  That $V$ on $[0,a]$ determines the asymptotics of $m$ up to the exponential error shown we'll call the \emph{direct direction} and the opposite we'll call the \emph{inverse direction}.  The full power of the $A$-function theory was used only for the inverse direction.
\end{remarks}

My main paper on this subject with Fritz~\cite{Fritz28} has many results including examples and a vast improvement in the proof of the direct direction of the local Borg-Marchenko theorem but I want to focus on what are undoubtedly, its two most important contributions.  First, unlike \cite{SiInv}, which is restricted to the finite $b$ case and to $V$ with strong conditions that imply $H$ is bounded from below, this paper allows \emph{any} $V$ in $L^1_{\text{loc}}$ including ones where $H$ is not bounded below.  It even allows the case where the operator is limit circle at $\infty$ with any choice of boundary condition at infinity.  The key observation is that one need only prove the direct direction of the local Borg-Markchenko theorem for $V$'s in that generality and then use the fact that one has proven the $A$-function basics for $L^1$ $V$'s of compact support.

The other issue addressed is how to determine $A(\alpha,x=0)$ from the spectral measure without reference to inverse Laplace transforms and \cite{Fritz28} proves the following formula
\begin{equation}\label{4.9}
  A(\alpha) = -2 \int_{-\infty}^{\infty} \lambda^{-\tfrac{1}{2}}\sin(2\alpha\sqrt{\lambda})\, d\rho(\lambda)
\end{equation}
where $\lambda^{-\tfrac{1}{2}}\sin(2\alpha\sqrt{\lambda})$ is defined for $\lambda\le 0$ by analytic continuation of the entire function from $\lambda>0$.  In one sense, this formula is poetry because the asymptotics of $m$ implies that $\int_{0}^{R}d\rho(\lambda)\sim\tfrac{2}{3\pi}R^{3/2}$ so the integral is never convergent.  However, \eqref{4.9} is proven in distributional sense and, in many cases, as an abelianized limit.

Fritz and I wrote a separate paper~\cite{Fritz29} proving both directions of the local Borg-Marchenko theorem without the full $A$-machinery and Benewitz~\cite{Bene} showed how one could modify and extend Borg's proof of his original theorem to obtain the local version.

%%%%%%%%%%%%%%%%%%%%%%%%%%%%%%%%%%%%%%%%%%%%%%%%%%%%%%%%%%%%%%
\section{Inverse spectral analysis with partial information on the potential} \lb{s5}
%%%%%%%%%%%%%%%%%%%%%%%%%%%%%%%%%%%%%%%%%%%%%%%%%%%%%%%%%%%%%%

My final section involves 5 papers with Fritz - three in a three part series~\cite{Fritz23, Fritz25, Fritz26} (one with del Rio) and two others~\cite{Fritz24, Fritz27}.  I can describe their theme by quoting two earlier theorems that motivated us (our work was just prior to 2000):

\begin{theorem} [Borg~\cite{Borg0}, 1946] \label{T5.1} Consider a Schr\"{o}dinger operator, $H(h_0,h_1)$, \eqref{1.3} on $L^2([0,1],dx)$ with $V\in L^1([0,1],dx)$ with boundary conditions
\begin{equation}\label{5.1}
  u'(0)+h_0u(0) = 0; \qquad u'(1)+h_1u(1) = 0
\end{equation}
with $h_0,h_1\in\bbR\cup\{0\}$.  Fix $h_1$ and two \emph{different} values, $h_0^{(a)}\ne h_0^{(b)}$ for $h_0$.  Then the eigenvalues of $H(h_0^{(a)},h_1)$ and $H(h_0^{(b)},h_1)$ determine $V$.
\end{theorem}

\begin{theorem} [Hochstadt-Lieberman~\cite{HochLieb}, 1978] \label{T5.2} With the notation of Theorem \ref{T5.1}, the spectra of one $H(h_0,h_1)$ and $V$ on $[0,\tfrac{1}{2}]$ determine $V$ on all of $[0,1]$.
\end{theorem}

All five papers explored this issue of what kind of combination of spectra and partial information on the parameters determine all the parameters.   \cite{Fritz24} discusses the Jacobi case and \cite{Fritz23} the Schr\"{o}dinger case on $\bbR$ with a.c. spectrum; in our discussion here, I'll focus on the other three papers which deal with the interval case discussed in those two earlier theorems.  Here are some examples in results in these three papers.

\begin{theorem} [Gesztesy-Simon~\cite{Fritz25}] \label{T5.3} Let $H(h_0,h_1)$ be as in Theorem \ref{T5.1}.  Then $h_0$, and $V$ on $[0,\tfrac{1}{2}+\tfrac{a}{2}]$ for some $a\in (0,1)$ and a set of eigenvalues, $S\subset\sigma(H)$ with
\begin{equation}\label{5.2}
  \#\{\lambda\in S\,\mid\,\lambda\le\lambda_0\} \ge (1-a) \#\{\lambda\in\sigma(H)\,\mid\,\lambda\le\lambda_0\}
\end{equation}
for all sufficiently large $\lambda_0$ uniquely determine $V$ on all of $[0,1]$ and determine $h_1$.
\end{theorem}

So colloquially speaking, an example is the potential on $[0,\tfrac{3}{4}]$ and half the eigenvalues determine $V$.

For the next result, we need some extra notation.  Given $c<d$, both in $\bbR$, we let $H(c,d;h_c,h_d)$ be the operator like $H(h_0,h_1)$ but on $L^2([c,d],dx)$ and with boundary conditions of \eqref{5.1} type at $c$ and $d$.

\begin{theorem} [Gesztesy-Simon~\cite{Fritz27}, Pivovarchik~\cite{Piv}] \label{T5.4} Let $0<a<1$ and let $S_-, S_+, S$ be three discrete infinite subsets of $\bbR$ bounded from  below which are pairwise disjoint. Let $h_a\in\bbR\cup\{\infty\}$ by given. If for some $L^1([0,1],dx)$ potential, V, and $h_0,h_1\in\bbR\cup\{\infty\}$, we have that $S_-, S_+, S$ are all the eigenvalues of $H(0,a;h_0,h_a), H(a,1;h_a,h_1), H(0,1;h_0,h_1)$ respectively, then $V$, $h_0$ and $h_1$ are uniquely determined.
\end{theorem}

\begin{remark} Pivovarchik~\cite{Piv} only has the case $a=\tfrac{1}{2}$.
\end{remark}

\begin{theorem} [del Rio-Gesztesy-Simon~\cite{Fritz26}] \label{T5.5} In the notation of Theorem \ref{T5.1}, $h_1$ is fixed, two-thirds of the spectra of $H(h_0,h_1)$ (in that the analog of \eqref{5.2} holds with $1-a=\tfrac{2}{3}$) for $3$ distinct values of $h_0$ determine $V$.
\end{theorem}

Without giving all the details, there are two ideas behind the proofs of these results.  First, if $m(z)$ is the $m$ function for solutions with $h_1$ BC at $1$ (i.e. given by \eqref{1.4} where $u_+$ is replaced by the solution obeying the BC at $1$), then $\lambda_0$ is an eigenvalue of $H(h_0,h_1)$ if and only if $m(\lambda_0)=-h_0$.  One knows about growth properties of $m(z)$ near infinity depending on the size of the interval $(b,c)$ on which we are considering the operator.  Results like Jensen's formula (see Simon~\cite[Section 9.8]{BCA}) imply that knowing the values of an entire function of a certain growth at a certain density of points determines the function.  Thus, eigenvalue information can determine $m$ which determines $V$ by Theorem \ref{T1.1}.  Second, one can use the knowledge of $m(\lambda_0,0)$, of $V(x)$ on $[0,a]$ and \eqref{1.7} to compute $m(\lambda_0,a)$.  Thereby, partial information on $V$ can transfer to information on an analytic function with less growth at infinity.

%%%%%%%%%%%%%%%%%%%%%%%%%%%%%%%

\end{document}